\documentclass[letterpaper,12pt]{article}  
\usepackage{amssymb}  
\usepackage{amsmath}  
\usepackage{makeidx}  
\usepackage{amscd}

\def\hfl#1#2#3{\smash{\mathop{\hbox  
to#3{\rightarrowfill}}\limits  
^{\scriptstyle#1}_{\scriptstyle#2}}}

\def\vfl#1#2#3{\llap{$\scriptstyle #1$}  
\left\downarrow\vbox to#3{}\right.\rlap{$\scriptstyle #2$}}

\newtheorem{theo}{Theorem}[section]    
\newtheorem{prop}[theo]{Proposition}  
\newtheorem{lem}[theo]{Lemma}  
\newtheorem{cor}[theo]{Corollary}    
\newtheorem{defi}[theo]{Definition}

\def \deg {{\rm{deg}}}

\def \Ga {{\Gamma}}  
\def \R {{\mathbb{R}}}  
\def \Pic {{\rm {Pic\,}}}

\def \Gal {{\rm{Gal\,}}}

\def \A{{\mathbb A}}  
\def \P{{\mathbb P}} 
  
\def \dim {{\rm{dim\,}}}

\def \GL {{\rm {GL}}}
\def \SL {{\rm {SL}}}  
\def \SO {{\rm {SO}}}

\def \Spin {{\rm {Spin}}}

\def\ov{\overline}

\def \Z {{\bf Z}}

\def \WW {{\rm W}}
\def \AA {{\rm A}}
\def \BB {{\rm B}}
\def \DD {{\rm D}}
\def \EE {{\rm E}}  
\def \RR {{\rm R}}  
\def \CC {{\rm C}}

\def \Bl {{\rm{Bl}}}

\def \id {{\rm{id}}}  
\def\G{{\bf G}}

\def\T{{\cal T}}

\def\H{{\rm H}}
  
\def\O{{\cal O}}

\def\Conv{{\rm Conv}}

\def\O{{\cal O}}

\def\Ga{\Gamma}

\def\St{{\rm St}}

\def\exp{{\rm exp}}
\def\g{{\mathfrak g}}
\def\p{{\mathfrak p}}
\def\n{{\mathfrak n}}
\def\h{{\mathfrak h}}

\newcommand{\bthe}{\begin{theo}}  
\newcommand{\ble}{\begin{lem}}  
\newcommand{\bpr}{\begin{prop}}  
\newcommand{\bco}{\begin{cor}}  
\newcommand{\bde}{\begin{defi}}  
\newcommand{\ethe}{\end{theo}}  
\newcommand{\ele}{\end{lem}}  
\newcommand{\epr}{\end{prop}}  
\newcommand{\eco}{\end{cor}}  
\newcommand{\ede}{\end{defi}}

\DeclareFontFamily{U}{wncy}{}  
\DeclareFontShape{U}{wncy}{m}{n}{%  
<5>wncyr5%  
<6>wncyr6%  
<7>wncyr7%  
<8>wncyr8%  
<9>wncyr9%  
<10>wncyr10%  
<11>wncyr10%  
<12>wncyr6%  
<14>wncyr7%  
<17>wncyr8%  
<20>wncyr10%  
<25>wncyr10}{}  
\DeclareMathAlphabet{\cyr}{U}{wncy}{m}{n}

\title{Del Pezzo surfaces and representation theory}  
\author{Vera Serganova and Alexei Skorobogatov} 
\date{\it{to Yuri Ivanovich Manin on his 70th birthday}} 
\begin{document}  
\baselineskip=15pt  

\maketitle  

\begin{abstract}

Manin discovered that del Pezzo surfaces 
are related to root systems. 
Batyrev conjectured that a universal torsor on a del Pezzo surface 
can be embedded into a certain projective homogeneous space of 
the semisimple group with the same root system, equivariantly with 
respect to the maximal torus action. 
Computational proofs of this conjecture based on the structure of 
the Cox ring have been recently given by Popov and Derenthal.
We give a new proof of Batyrev's conjecture 
using an inductive process interpreting
the blowing-up of a point on a del Pezzo surface in terms of
representations of 
Lie algebras corresponding to Hermitian symmetric pairs. 

\end{abstract}

\section*{Introduction}

Del Pezzo surfaces, classically defined as smooth surfaces of degree $d$
in the projective space $\P^d$, $d\geq 3$, 
are among the most studied and best understood 
of algebraic varieties. Over an algebraically closed ground field
such a surface is the quadric $\P^1\times \P^1$ or 
the projective plane $\P^2$ with $r=9-d$ points in
general position blown up; in this definition $d$ can be
any integer between 1 and 9. Despite the apparent simplicity
the enumerative geometry of these surfaces
displays amazing symmetries and puzzling coincidences. 
The 27 lines on a smooth cubic surface were discovered by Cayley and Salmon,
and the symmetries of their configuration were 
explored by Schoutte, Coxeter and Du Val. A breakthrough
came in the late 1960-s 
when Manin \cite{CF} discovered that to a del Pezzo surface $X$
of degree $d=9-r$, $d\leq 6$, one can attach 
a root system $\RR_r$ in such a way that 
the automorphism group of the incidence graph
of the exceptional curves on $X$ is the Weyl group $\WW(\RR_r)$. These
root systems are embedded into one another: $\RR_8=\EE_8$,
and as $r$ decreases one chops one by one the nodes off the long end of the
Dynkin diagram of $\EE_8$, until the diagram becomes disconnected. 
Let $\alpha_r$ be the simple root of $\RR_r$
corresponding to the node which must be removed 
from the Dynkin diagram of $\RR_r$ in order to obtain that
of $\RR_{r-1}$; let $\omega_r$ be the 
fundamental weight dual to $\alpha_r$.
For $r=4, 5, 6, 7$ the number of exceptional curves on $X$ is 
$|\WW(\RR_r)/\WW(\RR_{r-1})|=10, 16, 27, 56$,
respectively, and this is also
the dimension of the irreducible minuscule
representation $V(\omega_r)$ of the Lie algebra $\g_r$ of type $\RR_r$ 
with the highest weight $\omega_r$. It is tempting
to try to recover the Lie algebra directly from a del Pezzo
surface, but one has to bear in mind that the del Pezzo surfaces
of degree $d\leq 5$ depend on $10-2d$ moduli, so that the Lie
algebra should somehow take into account all del Pezzo surfaces
of given degree (see \cite{M}, and also \cite{FM}, \cite{Leung}).

Universal torsors were
introduced by Colliot-Th\'el\`ene and Sansuc in the 1970-s in a seemingly
unrelated line of research
(see \cite{CS} or \cite{S}). If $X$ is a smooth projective variety
over a field $k$,
then an $X$-torsor under a torus $T$ is a pair $(Y,f)$,
where $Y$ is a variety over $k$ with a free action of $T$, and $f$ is an
affine morphism $Y\to X$ whose fibres are orbits of $T$.
An $X$-torsor is universal if all invertible regular
functions on $Y$ are constant, and the Picard group of $Y$
is trivial (see Section \ref{S1} for details).
Then $T$ is isomorphic to
the N\'eron--Severi torus of $X$, i.e., the
algebraic torus dual to the Picard lattice of $X$
over an algebraic closure of $k$. 
In the work of Colliot-Th\'el\`ene, Sansuc,
Swinnerton-Dyer, Salberger and the second named author 
(see the references in \cite{S})
the birational geometry of universal torsors on del Pezzo surfaces 
of degrees 3 and 4 played
a crucial role in gaining some understanding of rational points on these
surfaces over number fields, for example, the Hasse principle, weak
approximation, the Brauer--Manin obstruction, the R-equivalence. 
The work of Batyrev, Tschinkel, Peyre, Salberger, Hassett, 
de la Bret\`eche, Heath-Brown, 
Browning and others on the Manin--Batyrev
conjectures on the number of rational points of bounded height,
highlighted the importance of explicitly describing 
universal torsors as algebraic varieties, and not merely
their birational structure. However, in the most interesting cases
such as those of (smooth) del Pezzo surfaces
of degrees 3 and 4, explicit equations of universal
torsors turned out to be quite complicated to write down explicitly.

Around 1990, Victor Batyrev told the second named author about his
conjecture relating universal torsors on del Pezzo surfaces
to certain projective homogeneous spaces.
Let $G_r$ be the simply connected semisimple group of type $\RR_r$.
We fix a maximal torus $H_r\subset G_r$, and a basis of simple roots
in the character group of $H_r$.
Let $P_r\subset G_r$ be the maximal parabolic subgroup defined by the root
$\alpha_r$ (the stabilizer of 
the line spanned by the highest weight vector of $V(\omega_r)$).
Batyrev conjectured that
a universal torsor $\T$ on a del Pezzo surface $X$ of degree $d=9-r$ 
over an algebraically closed field can be embedded into the affine cone 
$(G_r/P_r)_a\subset V(\omega_r)$ over $G_r/P_r$, equivariantly 
with
respect to the action of the N\'eron--Severi torus $T_r$ of $X$,
identified with an extension of $H_r$ by the scalar matrices $\G_m$. 
Moreover, the exceptional
curves on $X$ should be the images of the weight hyperplane
sections of $\T$ 
(that is, the intersections of $\T$ with the 
$H_r$-invariant hyperplanes in $V(\omega_r)$). 
Inspired by these ideas, the second named author
showed in \cite{S1} that the set of stable points of 
the affine cone over the Grassmannian $G(3,5)$
with respect to the action of the diagonal torus of $\SL(5)$, is a 
universal torsor over a del Pezzo surface of degree 5 which is
the GIT quotient by this action. 
Batyrev's line of attack 
on the general case of his conjecture uses the Cox ring
of $X$, which can be interpreted as the ring of regular
functions on a universal torsor over $X$. 
Indeed, Batyrev and Popov \cite{BP} (see also \cite{D1})
found the generators and relations of the Cox ring, which enabled Popov
in his unpublished thesis
\cite{P} in the case $d=4$ and Derenthal \cite{D} in the cases $d=3$ and
$d=2$ to prove Batyrev's conjecture by identifying
the generators with the weights of $V(\omega_r)$,
and comparing the relations with the well known
equations of $G_r/P_r$. The proofs of \cite{P} and \cite{D} 
are based on a substantial amount of calculation which  
grows rapidly with $r$, and do not seem to
give much insight into why things work this way.

In the present work we prove Batyrev's conjecture for del Pezzo surfaces
of degrees 4, 3 and 2 using a totally different approach, via
the representation theory of Lie algebras.
We start with the known case of a del Pezzo surface of degree 5.
(Alternatively, one could start with the simpler though somewhat irregular
case of degree 6, see \cite{BP}.) 
Let $\p_r$ be the Lie algebra of $P_r\subset G_r$.
We build an inductive
process based on the remark that the pair $(\RR_r,\alpha_r)$
for $r=4,5,6,7$ is a Hermitian symmetric pair, which is saying
that the complementary nilpotent algebra of $\p_r$ in 
$\g_r$ is commutative.
We show that $V(\omega_{r})$, as a $\g_{r-1}$-module,
has a direct factor isomorphic to $V(\omega_{r-1})$, and that 
the restriction of the projection $V(\omega_{r})\to V(\omega_{r-1})$
to a certain 
open subset $U\subset (G_r/P_r)_a$ is the composition of 
a $\G_m$-torsor and a morphism inverse to
the blowing-up of $V(\omega_{r-1})\setminus\{0\}$
at $(G_{r-1}/P_{r-1})_a\setminus\{0\}$ (see Corollary \ref{cc}).
Now we can explain the main idea of our proof.
Suppose that a universal torsor $\T$ over a
del Pezzo surface $X$ of degree $9-(r-1)$ is $T_{r-1}$-equivariantly embedded
into the affine cone $(G_{r-1}/P_{r-1})_a\subset
V(\omega_{r-1})$. Let $M$ be a point on $X$ outside of
the exceptional curves, and $\Bl_M(X)$ the blowing-up of $X$ at $M$. 
The space $V(\omega_{r-1})$ is a direct sum of $1$-dimensional
weight spaces of $H_{r-1}$, so that the torus consisting of
the diagonal matrices with respect to a weight basis of $V(\omega_{r-1})$
does not depend on the choice of this basis.
We show how to choose an element $t_M$ of this torus 
so that the translation $t_M^{-1}(G_{r-1}/P_{r-1})_a$ intersects $\T$ exactly
in the fibre of $f:\T\to X$ over $M$.
Then the closure of the inverse image 
of $t_M(\T\setminus f^{-1}(M))$ in $U$
is a universal torsor over $\Bl_M(X)$.
This yields a $T_r$-equivariant embedding of this universal torsor
into
$(G_{r}/P_{r})_a$. We then show that the image of this embedding
is contained in the open subset of stable points with a free
action of the N\'eron--Severi torus, so that $\Bl_M(X)$
embeds into the corresponding quotient.

Here is the structure of the paper. In Section \ref{S1}
we recall equivalent definitions and some basic properties of
universal torsors. In Section \ref{S2} we 
prove that the left action of a maximal torus of $G$ on $G/P$,
where $P$ is a maximal parabolic subgroup
of a semisimple algebraic group $G$, turns the set of stable points
with free action of the maximal torus into
a universal torsor on an open subset of the GIT quotient of 
$G/P$ by this action (with an
explicit list of exceptions, see Theorem \ref{univ}
for a precise statement).
In Section \ref{S3} we recall the necessary background from
the representation theory of Lie algebras.
The implications for the structure of the projection
of $(G_r/P_r)_a$ to $V(\omega_{r-1})$ are studied in Section \ref{S4}.
In Section \ref{DP} we list some well known properties of 
del Pezzo surfaces. Our main result, Theorem 
\ref{main}, is stated and proved in Section \ref{S5}.

The second named author is grateful to
the Centre de recherches math\'ematiques de l'Universit\'e 
de Montr\'eal, the Mathematical Sciences Research
Institute in Berkeley, and the organizers
of the special semester ``Rational and integral points on 
higher-dimensional varieties" for hospitality and support.

\section{Universal torsors}\label{S1}

Let $k$ be a field of characteristic $0$ with an algebraic closure $\ov k$.
Let $X$ be a geometrically integral variety over $k$.
We write $\ov X$ for $X\times_k\ov k$.
We denote by $\ov k[X]$ the $\ov k$-algebra of regular
functions on $\ov X$, and by $\ov k[X]^*$ the group of its invertible
elements.

Let $T$ be an algebraic $k$-torus, that is, an algebraic group
such that $\ov T\simeq \G_m^n$ for some $n$.
Let $\hat T\simeq\Z^n$ be the group of characters of $T$. 
The Galois group $\Ga=\Gal(\ov k/k)$ naturally acts on $\hat T$.

For generalities on torsors the reader is referred to \cite{S}.
An $X$-{\sl torsor} under $T$ is a pair $(\T,f)$, where $\T$
is a $k$-variety with an action of $T$,
and $f:\T\to X$ is a morphism such that
locally in \'etale topology $\T\to X$ is $T$-equivariantly
isomorphic to $X\times_k T$.
The following lemma is well known.

\ble \label{torsor}
Suppose that a $k$-torus $T$ acts on a $k$-variety $Y$ with trivial stabilizers,
and $g:Y\to X$ is an affine morphism of $k$-varieties whose
fibres are orbits of $T$. Then $g:Y\to X$ is a torsor under $T$.
\ele
{\em Proof}
The property of $g$ to be a torsor can be checked locally on $X$. 
Let $U$ be an open affine subset of $X$.
Since $g$ is affine, $g^{-1}(U)$ is also affine
(\cite{H}, II, 5, Exercise 5.17). Since the stabilizers of 
all $\ov k$-points of $g^{-1}(U)$ are trivial, by a corollary of
Luna's \'etale slice theorem (see \cite{GIT}, p. 153)
the natural map $g^{-1}(U)\to U$ is a torsor under $T$.
The lemma follows. QED

\medskip

Colliot-Th\'el\`ene and Sansuc associated to a torsor $f:\T\to X$
under a torus $T$ the exact sequence
of $\Ga$-modules (\cite{CS}, 2.1.1)
\begin{equation}
1\to \ov k[X]^*/\ov k^*\to \ov k[\T]^*/\ov k^*\to \hat T\to
\Pic\ov X\to\Pic\ov \T\to 0.\label{cts}
\end{equation}
Here the second and the fifth arrows are induced by $f$.
The forth arrow is called the {\sl type} of $\T\to X$. To define it
consider the natural pairing compatible with the action of 
the Galois group $\Ga$:
$$\cup:\ \H^1(\ov X,T)\times \hat T\to \H^1(\ov X,\G_m)=\Pic \ov X,$$
where the cohomology groups are in \'etale or Zariski
topology. The type sends $\chi\in \hat T$ to $[\ov \T]\cup\chi$,
where $[\ov \T]\in \H^1(\ov X,T)$ is the class of the torsor
$\ov \T\to \ov X$. A torsor $\T\to X$ is called {\sl universal}
if its type is an isomorphism. If the variety $X$ is projective,
(\ref{cts}) gives the following characterisation of universal
torsors: an $X$-torsor under a torus is universal if and only if 
$\Pic\ov \T=0$ and $\ov k[\T]^*=\ov k^*$, that is, $\ov \T$ has no
non-constant invertible regular functions.

We now give an equivalent definition of type which
does not involve cohomology. Let $K=\ov k(X)$ be the function 
field of $\ov X$, and $\T_K$ the generic fibre of $\ov \T\to \ov X$.
By Hilbert's Theorem 90 the $K$-torsor $\T_K$ is trivial, that is, is
isomorphic to $T_K=T\times_kK$. By Rosenlicht's lemma we have
an isomorphism of $\Ga$-modules $K[\T_K]^*/K^*=K[T_K]^*/K^*=\hat T$.
This isomorphism associates to a character $\chi\in \hat T$ a rational function
$\phi\in \ov k(\T)^*$ such that $\phi(tx)=\chi(t)\phi(x)$;
the function $\phi$ is well defined up to
an element of $K^*=\ov k(X)^*$. The divisor of $\phi$ on $\ov \T$
does not meet the generic fibre $\T_K$, and hence comes from
a divisor on $\ov X$ defined up to a principal divisor.
We obtain a well defined class $\tau(\chi)$ in $\Pic \ov X$.

\ble \label{tau}
The map $\tau:\hat T \to \Pic \ov X$ coincides with the type of
$f:\T\to X$ up to sign.
\ele
{\em Proof} According to \cite{S}, Lemma 2.3.1 (ii),
the type associates to $\chi$ the subsheaf $\O_\chi$
of $\chi$-semiinvariants of the sheaf $f_*(\O_\T)$. 
The function $\phi$ is a rational
section of $\O_\chi$, hence the class of its divisor 
represents $\O_\chi\in \Pic \ov X$. QED

\medskip

For the sake of completeness we note that if $f:\T\to X$ 
is a universal torsor, then 
the group of divisors on $\ov X$ is naturally identified with
$K[\T_K]^*/\ov k^*$; this identifies the semigroup of effective divisors
on $\ov X$ with $(K[\T_K]^*\cap \ov k[\T])/\ov k^*$.

We have $\ov k[\T]=\oplus_{\chi\in\hat T}\ov k[\T]_\chi$,
where $\ov k[\T]_\chi$ is the set of regular functions $\phi$
on $\T$ satisfying $\phi(tx)=\chi(t)\phi(x)$ for any $t$ in $T$.
We also define $\ov k(\T)_\chi$ as the set of rational functions
on $\T$ satisfying the same condition.
Since $\ov k(\T)_\chi$ is the group of rational sections of 
the sheaf $\O_\chi$, we have $\ov k[\T]_\chi=\H^0(\ov X,\O_\chi)$.
Hence if the sheaf $\O_\chi$ defines a morphism 
$X\to\P(\H^0(\ov X,\O_\chi)^*)$, we obtain a commutative
diagram
\begin{equation}
\begin{array}{cclc}\T&\to&\ov k[\T]_\chi^*\setminus\{0\}=&
\H^0(\ov X,\O_\chi)^*\setminus\{0\}\\
\downarrow&&&\downarrow\\
X&&\to&\P(\H^0(\ov X,\O_\chi)^*)
\end{array}\label{end}
\end{equation}

\section{$G/P$ and the torus quotient}\label{S2}

Let $G$ be a split simple simply connected algebraic group over $k$,
with a split maximal torus $H\subset G$; in this case the root system
$\RR$ of $G$ relative to $H$ is irreducible.
Write $\hat H$ for the character group of $H$.
We use the standard notation $Q(\RR)$
for the lattice generated by the simple roots, then $P(\RR)=\hat H$ is the
dual lattice generated by the fundamental weights. We denote the Weyl group
by $\WW =\WW (\RR)$. 

Let $G\to\GL(V)$ be an irreducible 
representation of $G$ with a fundamental highest weight $\omega\in \hat H$.
Let $v\in V$ be a highest weight vector.  
The stabilizer of the line $kv$ is a maximal parabolic subgroup
$P\subset G$. 
The homogeneous space $G/P$ is thus a smooth projective subvariety 
of $\P(V)$ (indeed, the
only closed orbit of $G$ in $\P(V)$). We write $\hat P$ (resp.
$\hat G$) for the character group of $P$ (resp. of $G$).
Let $\varepsilon:\hat P\to \Pic G/P$ be the map
associating to the character
$\chi\in \hat P$ the $G/P$-torsor under $\G_m$ defined
as the quotient of $G\times \G_m$ by $P$,
where $p\in P$ sends $(g,t)$ to $(gp^{-1},\chi(p)t)$.
This map fits into the exact sequence
$$0\to \hat G\to \hat P\to \Pic G/P\to \Pic G\to 0.$$
Since $G$ is semisimple and simply connected
we have $\hat G=\Pic G=0$, so that
$\varepsilon$ is an isomorphism (see, e.g., \cite{VP}). 
Since $\hat P$ is the subgroup of 
$\hat H$ generated by $\omega$, we see that $\Pic G/P$
is generated by the hyperplane section class. This fact implies
the following elementary statement from
projective geometry.

\ble \label{codim2}
Let $L_1$ and $L_2$ be distinct hyperplanes in the projective space $\P(V)$.
Then $(G/P)\cap L_1\cap L_2$ has codimension $2$ in $G/P$.
\ele
{\em Proof} Since $\Pic(G/P)$ 
is generated by the class of a hyperplane section,
for any hyperplane $L\subset \P(V)$ the closed subset $(G/P)\cap L$
is irreducible of codimension 1, and the intersection has multiplicity 1.
If the codimension of $(G/P)\cap L_1\cap L_2$ in $G/P$ is 1, we
have $(G/P)\cap L_1\cap L_2=(G/P)\cap L$
for any $L$ in the linear family spanned by $L_1$ and $L_2$.
Choosing $L$ passing through a point of $G/P$ not contained in $L_1$,
we deduce a contradiction. QED

\medskip

By the irreducibility of $V$ the centre
$Z(G)$ acts diagonally on $V$, and hence it acts trivially on $\P(V)$.
For a $\ov k$-point $x\in \P(V)$
we denote the stabilizer of $x$ in $H$ by $\St_H(x)$.
We now show that
for $x$ in a dense open subset of $G/P$ we have 
$\St_H(x)=Z(G)$, and determine the points 
such that $\St_H(x)$ is strictly bigger than $Z(G)$.

\bpr \label{free}
Let $x$ be a $\ov k$-point of $G/P$, and let $K_x$ be the 
connected component of the centralizer of $\St_H(x)$ in $G$. Then
we have the following properties:

{\rm (i)} $K_x$ is a reductive subgroup of $G$, $H\subset K_x$;

{\rm (ii)} $x\in K_x wv=K_x/(wPw^{-1}\cap K_x)$ for some $w\in \WW$;

{\rm (iii)} $Z(K_x)=\St_H(x)$;

{\rm (iv)} $\St_H(x)$ is finite if and only if $K_x$ is semisimple, in 
which case the ranks of $K_x$ and $G$ are equal.
\epr
{\em Proof} If $\St_H(x)=Z(G)$, then $K_x=G$, and all the statements are
clearly true. Assume that $\St_H(x)$ is bigger than $Z(G)$, then
$K_x$ is a closed subgroup of $G$, $K_x\not=G$.

Let ${\mathfrak k}_x$ be the Lie algebra of $K_x$; explicitly
${\mathfrak k}_x\subset\g$ is the fixed set of
${\rm Ad}(\St_H(x))$. Since ${\mathfrak k}_x$ contains
the Cartan subalgebra ${\mathfrak h}$, it has a root
decomposition ${\mathfrak k}_x=\oplus_{\alpha\in S}\g_\alpha$,
where $S\subset\RR$. Let $\exp_\alpha\in \hat H $ be the 
multiplicative character defined by the root $\alpha\in\RR$.
The space $\g_\alpha$ consists of $y\in\g$
such that ${\rm Ad}(h)y=\exp_\alpha(h)y$ for all $h\in H$. Thus 
$\g_\alpha\subset{\mathfrak k}_x$ if and only if 
$\St_H(x)\subset H$ is in the kernel of $\exp_\alpha$.
Therefore $S=-S$, so that ${\mathfrak k}_x$ is reductive, and
hence so is $K_x$.

The fixed points of $H$ in $G/P$ come from the points
$wv$, where $w\in\WW $.
One of these, say $x_0=wv$, is contained in the closure of the orbit
$Hx$. The stabilizer of $x_0$ in $G$ is the parabolic subgroup
$wPw^{-1}$. To prove (ii) we 
need to show that $x$ belongs to the $K_x$-orbit of $x_0$.
Let $N\subset G$ be the unipotent subgroup complementary
to $wPw^{-1}$, that is, such that the corresponding Lie
algebras satisfy $\g=\n\oplus w\p w^{-1}$.
Then $N\cap wPw^{-1}=\{1\}$, and the $N$-orbit of the line $kx_0$
is the open Schubert cell $Nx_0\subset G/wPw^{-1}\simeq G/P$.
The intersection of this open Schubert cell with $Hx$
is a non-empty open subset of $Hx$, thus there is 
a $\ov k$-point $x_1\in Hx\cap Nx_0$.
We can write $x_1=u.x_0$ for some $u\in N$.
The complement to the union of connected components
of the centralizer of $\St_H(x)$ other than $K_x$, is an
open neighbourhood of $1$ in $G$. We choose $x_1$ 
in such a way that $u$ belongs to this open set.
Since $H\subset K_x$, the points $x$ and $x_1$ are in the same $K_x$-orbit,
so that it is enough to show that $x_1\in K_x x_0$.
Any $t\in \St_H(x)$ fixes both $x_1$ and $x_0$, thus
$x_1=u.x_0=t^{-1}ut.x_0$. Therefore, $u^{-1}t^{-1}ut$ fixes $x_0$,
hence $u^{-1}t^{-1}ut\in wPw^{-1}$. On the other hand, $H$ normalizes $N$,
thus $t^{-1}ut\in N$, implying $u^{-1}t^{-1}ut\in N$. Since the
intersection of $wPw^{-1}$ and $N$ is $\{1\}$, we see that $u$ and $t$
commute. By the choice of $x_1$ we see that $u$ is in the connected component
of $1$ of the centralizer of $\St_H(x)$, that is,
$u\in K_x$. This completes the proof of (ii).

The centre of $K_x$ is contained in every maximal torus, in
particular, in $H$. Any element of $Z(K_x)$
fixes $x$, since $x\in K_x/(wPw^{-1}\cap K_x)$, so that $Z(K_x)\subset\St_H(x)$.
On the other hand, every element of $\St_H(x)$ commutes with $K_x$
by the definition of $K_x$. But $\St_H(x)\subset H\subset K_x$, hence
$\St_H(x)\subset Z(K_x)$. This proves (iii).

The rank of the semisimple part of $K_x$ equals the rank of $G$
if and only if $Z(K_x)$ is finite. If $Z(K_x)$ is finite, 
then $K_x$ is semisimple by definition. Thus (iv) follows from (iii).
QED
\medskip

Let us fix a weight basis in $V$, that is, a basis in which $H$
is diagonal. The weight of a coordinate is the character
of $H$ by which $H$ acts on it. Denote by $\Lambda$ the set of 
weights of $H$ in $V$, and by ${\rm wt}(x)$ 
the set of weights of $x\in G/P$,
that is, the weights of the non-vanishing coordinates of $x$.

\bco \label{c1}
Assume that $\RR$ is simply laced.
Then the codimension of the set of $\ov k$-points $x\in G/P$ 
such that $\St_H(x)$ is finite, and $\St_H(x)\not=Z(G)$, 
is at least $2$.
\eco
{\em Proof} By Proposition \ref{free} and $\WW$-invariance it is 
sufficient to show that the codimension of $Kv$ in $Gv$ is at least
2 for any proper connected semisimple subgroup $K\subset G$ containing $H$.
(The set of such subgroups is clearly finite.)

For any $x\in G/P$ the property ${\rm wt}(x)=\Lambda$ implies $\St_H(x)=Z(G)$.
Let $V'\subset V$ be the irreducible representation of $K$ generated by $v$.
Denote by $\Lambda'$ the set of weights of $V'$, and write $V=V'\oplus U$,
where $U$ is another $K$-invariant subspace. First, we claim that
$\Lambda'\neq\Lambda$ because otherwise one can find 
$x\in \mathbb{P}(Kv)$ 
such that ${\rm wt}(x)=\Lambda$, and $\St_H(x)=Z(G)=Z(K)$ would imply 
$K=G$. In particular, $U\neq 0$. If $\dim U>1$, then
the codimension of $Kv \subset Gv\cap V'$ is at least 2 by Lemma ~\ref{codim2}.

If $\dim U=1$, then $U$ is a trivial representation
of $K$ and $0$ is not a weight of $V'$. But then $U$ is invariant 
under the action of the Weyl group $\WW$. Therefore $wKw^{-1}$ acts trivially
on $U$ for any $w\in \WW$. If $a\in \RR$ is a root of $K$, then $w(a)$ 
is a root of $wKw^{-1}$, but in the simply laced case $\WW$ acts 
transitively on $\RR$, hence the subgroups $wKw^{-1}$, $w\in \WW$,
generate the whole group $G$. Thus, $U$ is $G$-invariant, but that contradicts
the irreducibility of $V$. QED

\medskip

Recall that a $\ov k$-point $x\in V$ is called {\sl stable}
for the action of $H$ if the orbit $Hx$ is closed,
and the stabilizer of $x$ in $H$ is finite (\cite{GIT}, p. 194).
We always consider stability with respect to the action of $H$,
and drop the reference to $H$ when it causes no confusion.

For a subset $M\subset\hat H$ we write $\Conv(M)$ for the convex
hull of $M$ in the vector space $\hat H\otimes \R$.
It is well known that $\Conv(\Lambda)=\Conv(\WW\omega)$
(\cite{GS}, \cite{FH}, see \cite{Da}, Prop. 2.2 (i) for a short proof).
The Hilbert--Mumford
numerical criterion of stability says that $x$ is stable 
if and only if $0$ belongs to the interior of 
$\Conv({\rm wt}(x))$  (\cite{Do}, Thm. 9.2). 

In the following statement and thereafter the numeration of the nodes
of Dynkin diagrams, simple roots and
fundamental weights follows the conventions of \cite{B}.

\bpr \label{stable}
Assume that the pair $(\RR,\omega)$ is not in the following list:
\begin{equation}
(\RR_r,\omega_1),\ (\AA_r,\omega_r), \ (\AA_3,\omega_2), \ 
(\BB_2,\omega_2),\ (\CC_2,\omega_2),\
(\DD_4,\omega_3), \ (\DD_4,\omega_4),\label{list}
\end{equation}
where $\RR_r$ is $\AA_r$, $\BB_r$, $\CC_r$, or $\DD_r$.
Let $x$ be a point of $V\otimes_k\ov k$ such that no two elements of
$\WW\omega\setminus {\rm wt}(x)$ are adjacent vertices of 
$\Conv(\WW\omega)$. Then $x$ is stable. 

In particular, the set of unstable points of $G/P$
has codimension at least $2$.
\epr
{\em Proof} Since $\sum_{w\in\WW}w\omega=0$, the point $0$
is contained in the interior of $\Conv(\WW\omega)=\Conv(\Lambda)$
in $\hat H\otimes \R$. Thus
if all the coordinates of $x$ with weights in $\WW\omega$
are non-zero, then $x$ is stable. 
  
Now assume that exactly one such coordinate of $x$ is zero; 
because of the action of $\WW$ it is no
loss of generality to assume that it corresponds to
$\omega$. The dimension of the corresponding eigenspace is 1,
so to check that $x$ is stable it is enough to show that
$0$ lies in the interior of $\Conv(\WW \omega\setminus\{\omega\})$.
The vertices of $\Conv(\WW \omega)$
adjacent to $\omega$ are $\omega-w\alpha$,
where $\alpha$ is the root dual to $\omega$,
for all $w$ in the stabilizer of $\omega$ in $\WW $
(see \cite{FH}, Lemma 3 and Cor. 2). All these are contained
in the hyperplane $L=0$, where
$$L(y)=(y,\omega)-(\omega^2)+(\omega,\alpha)=
(y,\omega)-(\omega^2)+\frac{1}{2}(\alpha^2).$$
We have $L(\omega)>0$. Thus $0$ belongs to the interior of 
$\Conv(\WW \omega\setminus\{\omega\})$ if and only if $\omega$ and $0$
are separated by this hyperplane, that is, if and only if
$L(0)<0$. Therefore, we need to check the condition
$$
(\omega^2)>\frac{1}{2}(\alpha^2).
$$
Note that the numbers $2(\omega^2)/(\alpha^2)$, 
for all possible fundamental weights,
are the diagonal elements of the inverse Cartan matrix of $\RR$.
A routine verification using the tables of \cite{B} or \cite{VO}
shows that this inequality is satisfied for the pairs $(\RR,\omega)$
not in the list (\ref{list}).

Finally, let $\WW \omega\setminus{\rm wt}(x)=\{\lambda_1,\ldots,\lambda_n\}$.
By assumption $\lambda_1,\ldots,\lambda_n$ correspond
to pairwise non-adjacent
vertices of $\Conv(\WW \omega)$.
Thus 
$$\Conv(\WW \omega\setminus \{\lambda_1,\ldots,\lambda_n\})=
\bigcap_{i=1}^n \Conv(\WW \omega\setminus \{\lambda_i\}).$$
Since $0$ is in the interior of each convex hull in the
right hand side, it is also in the interior of $\Conv({\rm wt}(x))$. 

The last statement is an application of Lemma \ref{codim2}. QED

\bde
Let $T\subset \GL(V)$ be the torus 
generated by the image of $H$ in $\GL(V)$ and the scalar matrices 
$\G_m\subset \GL(V)$.
We write $(G/P)_{a}$ for the affine cone over $G/P$ in $V$, and
$(G/P)_{a}^{sf}$ for the open subset of stable points with trivial stabilizers
in $T$. 
\ede 
By the irreducibility of $V$,
the stabilizer of $x\in V\otimes_k\ov k$, $v\not=0$, in $T$ is trivial
if and only if
$\St_H(pr(x))=Z(G)$, where $pr(x)$ is the image of $x$ in $\P(V)$.

\ble \label{Y}
There exist a smooth quasi-projective variety $Y$ and an affine
morphism $f : (G/P)_{a}^{sf}\to Y$ which is a torsor with 
structure group $T$ with respect to its natural left action on $G/P$. 
\ele
{\em Proof} 
By the geometric invariant theory there
exist a quasi-projective variety $Y$ and an affine morphism 
$f : (G/P)_{a}^{sf}\to Y$ such that every fibre of $f$ is an orbit of $T$
(\cite{GIT}, Thm. 1.10 (iii)).
Since the stabilizers of all $\ov k$-points of $(G/P)_{a}^{sf}$
are trivial, Lemma \ref{torsor} implies that
$f : (G/P)_{a}^{sf}\to Y$ is a torsor under $T$. 
The smoothness of $Y$ follows from the smoothness of $(G/P)_{a}$,
since a torsor is locally trivial in \'etale topology. QED

\bthe \label{univ}
Assume that the root system $\RR$ is simply laced,
and the pair $(\RR,\omega)$
is not in the list $(\ref{list})$.
Then the only invertible regular
functions on $(G/P)_{a}^{sf}$ are constants, so that 
$f : (G/P)_{a}^{sf}\to Y$ is a universal torsor.
\ethe
{\em Proof} By Lemma \ref{Y} we need to show that $\Pic \ov \T=0$
and $\ov k[\T]^*=\ov k^*$ for $\T=(G/P)_{a}^{sf}$ (see Section \ref{S1}).
The Picard group of $(G/P)_{a}$
is trivial since that of $G/P$ is generated by the 
class of a hyperplane section. Thus it suffices to show that the complement
to $(G/P)_{a}^{sf}$ in $(G/P)_{a}$ has codimension at least 2.
The set of unstable points has codimension at least 2, by
Proposition \ref{stable}. The closed subset of its complement
consisting of the stable points
with non-trivial (finite) stabilizers in $T$, also has
codimension at least 2, as follows from Corollary \ref{c1}.
QED

\section{Hermitian symmetric pairs}\label{S3}

Let $\g$ be a semisimple Lie algebra over the field $k$
with Chevalley basis
$\{H_\beta,X_\gamma\}$, where $\gamma$ is a root of $\RR$, 
and $H_\beta=[X_\beta,X_{-\beta}]$, 
where $\beta$ is a simple root of $\RR$. 

A simple root $\alpha$ of $\g$ defines
a $\Z$-grading on $\g$ in the following way.
We set $\deg(X_\alpha)=1$, $\deg(X_{-\alpha})=-1$,
$\deg(X_{\pm \beta})=0$ for all other simple roots $\beta\not=\alpha$,
and $\deg(H_\beta)=0$ for
all simple roots $\beta$. Then 
\begin{equation}
\g=\bigoplus_{i=-l(\alpha)}^{l(\alpha)}\,\g_i, \label{gr}
\end{equation} 
where
$l(\alpha)$ is the {\sl label} of $\alpha$, that is, the coefficient of $\alpha$
in the decomposition of the maximal root as a linear combination
of the simple roots. The Lie algebra $\p=\oplus_{i\geq 0}\g_i$
is the parabolic subalgebra defined by $\alpha$, and 
$\n=\oplus_{i<0}\g_i$ is the complementary
nilpotent algebra. The centre of the
Lie algebra $\g_0$ is one-dimensional, so that $\g_0=Z(\g_0)\oplus \g'$,
where $\g'$ is the semisimple Lie algebra whose Dynkin diagram
is that of $\g$ with the node corresponding to $\alpha$ removed.

It is clear from (\ref{gr}) that $l(\alpha)=1$ if and only if $[\n,\n]=0$.
The following terminology has its origin in the theory of symmetric spaces,
see \cite{He}, Ch. VIII.

\bde
The pair $(\RR,\alpha)$ is a Hermitian symmetric pair if $l(\alpha)=1$,
or, equivalently, if $\n$ is a commutative Lie algebra.
\ede

If $\RR$ is simply laced, then $(\RR,\alpha)$
is a Hermitian symmetric pair if and only if $\RR=\AA_n$, or
if it is one of the following pairs:
$(\DD_n,\alpha_i)$, where $i=1$, $n-1$ or $n$, $(\EE_6,\alpha_1)$,
$(\EE_6,\alpha_6)$, and $(\EE_7,\alpha_7)$.

\medskip

We now assume that $\n$ is commutative.
Our next goal is to explore the implications of this assumption
for the restriction of the $\g$-module $V$
to the semisimple subalgebra $\g'$. We write
${\rm U}({\mathfrak l})$ for the universal enveloping
algebra of the Lie algebra ${\mathfrak l}$, and
$S(W)$ for the symmetric algebra of the vector space $W$.
Since $\n$ is commutative we have ${\rm U}(\n)=S(\n)$.

The line $kv$ is a 1-dimensional $\p$-submodule of $V$, hence
the $\g$-module $V$ is the quotient of the induced module
${\rm U}(\g)\otimes_{{\rm U}(\p)}kv$
by the submodule generated by $X_{-\alpha}^2v$. 
(This follows from the construction of $V$ as the quotient of the
Verma module by the submodule generated by $X_{-\beta}v$ for the simple roots
$\beta\not=\alpha$, and $X_{-\alpha}^2v$.)
By the Poincar\'e--Birkhoff--Witt theorem we have
${\rm U}(\g)={\rm U}(\p)\otimes_k{\rm U}(\n)$.
The line $kv$ is a trivial $\g'$-module.
Therefore, the $\g'$-module ${\rm U}(\g)\otimes_{{\rm U}(\p)}kv$
is isomorphic to ${\rm U}(\n)=S(\n)$, so that
the finite dimensional vector space
$V$ inherits the $\Z_{\leq 0}$-graded commutative $k$-algebra structure from
$S(\n)$, $V=\oplus_{n\leq 0}V^n$.
We turn this grading into a $\Z_{\geq 0}$-grading by
setting $V_n=V^{-n}$.
Since $\g'$ has grading $0$, the direct sum
$V=\oplus_{n\geq 0}V_n$ is the direct sum
of $\g'$-modules, and we can write
$$V=k\oplus\n\oplus \, S^{\geq 2}(\n)/S(\n){\rm U}(\g')
X_{-\alpha}^2\,,$$
where $k=V_0$, $\n=V_1$. Note that $1\in V_0$ is a highest
weight vector; it generates $V$ as a $S(\n)$-module.

\ble \label{rep}
Let $(\RR,\alpha)$ be a Hermitian symmetric pair. 
Then the adjoint representation of $\g'$ on $V_1=\n=\g_{-1}$
is the irreducible representation such that $X_{-\alpha}$
is a highest weight vector. If $\RR$
is simply laced, then the highest weight
$\omega'$ of $V_1$ is the sum of the fundamental
weights corresponding to the nodes of the Dynkin diagram of 
$\RR$ adjacent to the node $\alpha$.
\ele
{\em Proof} We have
$[X_{\beta},X_{-\alpha}]=0$ for all simple roots $\beta\not=\alpha$, so that
$X_{-\alpha}$ is annihilated by the positive roots of 
$\g'$. Every root of $\n$ is the sum of
$-\alpha$ and a root of $\g'$, so that $\n$ is generated
by $X_{-\alpha}$ as a $\g'$-module. The computation of the
weight of $X_{-\alpha}$ is immediate from the defining relations
among the elements of the Chevalley basis. QED

\medskip

We have the exponential map 
$$\exp : \n\to S(\n), \quad \exp(u)=1+u+\frac{1}{2}u^2+\frac{1}{3!}u^3+\ldots$$
Let $G$ be the simply connected semisimple algebraic $k$-group
with Lie algebra $\g$, $P\subset G$ the parabolic subgroup with Lie
algebra $\p$, and $N$ the unipotent $k$-group with Lie algebra $\n$.
By the Chevalley construction of the Lie group from its Lie
algebra, $N$ acts on $V$ by the rule $1+x\mapsto \exp(x)$.
Recall that the open Schubert cell of $G/P\subset \P(V)$
is the $N$-orbit of the highest weight vector, and hence
is identified with $\exp(\n)$. 
(In particular, $\dim G/P=\dim V_1$.)
Thus $\exp(x)$ is a polynomial $G'$-equivariant map
$$\exp:V_1\to (G/P)_a\subset V=\oplus_{n\geq 0}V_n.$$
Let $p : V_1=\n\to V_2$ 
be the degree $2$ graded component of $\exp(x)$.

\ble \label{p}
Let $G'$ be the simply connected semisimple $k$-group with the Lie
algebra $\g'$, and $P'\subset G'$ the parabolic subgroup
which is the stabilizer of the line spanned by the highest
weight vector $X_{-\alpha}\in\n$. 
The restriction of $\exp(x)$ to $(G'/P')_a$ coincides with
$(1,\id,0,0,\ldots)$. We have $(G'/P')_a=p^{-1}(0)$, and
the ideal of $(G'/P')_a$ is generated by the coordinates of $p(x)$.
\ele
{\em Proof} It is clear that every graded component of $\exp(x)$
of degree at least 2 sends the orbit $(G'/P')_a$ of the highest
weight vector $X_{-\alpha}$ to $0$. Indeed, 
$X_{-\alpha}^m$ is in the 
kernel of the natural map $S^m(\n)\to V_m$, for $m\geq 2$ . 
To prove the second statement let us observe that
the symmetric square $S^2(\n)$ decomposes as
the direct sum of $V_2$ and the $\g'$-submodule generated 
by $X_{-\alpha}^2$, which is the irreducible representation $V(2\omega')$
with
highest weight $2\omega'$. It is well known (\cite{LT}, proof of Thm. 1.1,
or \cite{BP}, Prop. 4.2) that the orbit of the highest
weight vector is the intersection of the second Veronese embedding with
$V(2\omega')$. This completes the proof. QED

\medskip

Consider the following series of root systems:
\begin{equation}
\AA_{4}\subset \DD_{5}\subset \EE_{6}\subset \EE_{7}. \label{series}
\end{equation}
Let $(\RR,\alpha)$ be one of the Hermitian symmetric pairs
\begin{equation}
(\AA_4,\alpha_3), \quad (\DD_5,\alpha_5), \quad (\EE_6,\alpha_6), \quad 
(\EE_7,\alpha_7), \label{herm}
\end{equation} 
where the roots are numbered as in \cite{B}. By Lemma \ref{rep}
the pair $(G',P')$ is defined by $(\RR',\alpha')$
which is the previous pair to $(\RR,\alpha)$ in $(\ref{herm})$.
In other words, $P'$
corresponds to the only node of the smaller diagram adjacent to
$\alpha$. (If $G$ is of type $\AA_4$, then 
$G'$ is of type $\AA_1\times\AA_2$, $G'/P'\simeq\P^1\times\P^2$, but 
we shall not have to consider this case.) 

We note that the fundamental weight $\omega$ dual to $\alpha$ 
is minuscule, that is, the 
weights of $V$ are $\WW \omega$, and $\WW v$ is a basis
of $V$ (see \cite{B}, VIII.7.3). 
We also note that the $G$-module $V$
defined by $\omega$ is faithful (this follows from the fact that
$\omega$ generates $P(\RR)/Q(\RR)$, which can be checked from the tables).
Thus the faithful
representation of $G$ in $V$ defines a faithful representation of 
$G'$, and this implies that $G'\subset G$ (in fact, $G'$ is 
the Levi subgroup of $P$). 

\medskip

Let us identify the graded components of $V$ in various cases.
Let $d_r=\dim V$. We have
$$d_4=10, \ d_5=16, \ d_6=27, \ d_7=56.$$
The details given below show that
for $r=4,5,6$ the graded components of $\exp(x)$
of degree at least 3 are zero.

Let $\RR=\AA_4$. Then $G=\SL(5)$, and $G/P$ is the Grassmannian
$G(2,5)$. Let us denote by $E_n$ the standard $n$-dimensional 
representation of $\SL(n)$. We have
$V=\Lambda^2(E_5)$, $\dim V=10=1+6+3$. 
The group $G'=\SL(2)\times\SL(3)$ is
embedded into $\SL(5)$ in the obvious way, and
the graded factors of $V$ are $V_1=E_2\otimes E_3$, 
$V_2=\Lambda^2(E_3)\cong E_3^*$. 
The map $p : V_1\to V_2$ sends $x$ to the 
$\Lambda^2(E_3)$-component of $$x\wedge x\in \Lambda^2(E_5)=\Lambda^2(E_2)
\oplus (E_2\otimes E_3)\oplus\Lambda^2(E_3).$$

Let $\RR=\DD_5$. Then $V$ is a spinor representation of
$G=\Spin(10)$ of dimension $16=1+10+5$, and $G/P$ is the isotropic Grassmannian
(one of two families of maximal isotropic subspaces
of the non-degenerate quadric of rank 10), $\dim G/P=10$.
The graded components are $V_1=\Lambda^2(E_5)$ and
$V_2=\Lambda^{4}(E_5)\cong E_5^*$.
The map $p : V_1\to V_2$ sends $x$ to $x\wedge x$.

Let $\RR=\EE_6$. Then $\dim V=27=1+16+10$, $V_1$ is the spinor
representation of $\Spin(10)$ as above, and $V_2$
is the standard $10$-dimensional representation of $\SO(10)$. 
We have $\dim G/P=16$.

Let $\RR=\EE_7$. Then $\dim V=56=1+27+27+1$, $V_1$
is the $27$-dimensional representation of the group
of type $\EE_6$ considered above, $V_2=(V_1)^*$, 
and $V_3=k$ is the trivial 1-dimensional representation.
(The graded components of degree at least 4 are zero.)
We have $\dim G/P=27$. We define
$q : V_1=\n\to V_3=k$ as the degree $3$ graded component of ${\rm exp}(x)$.
This is a $\EE_6$-invariant cubic form in 27 variables.
The 27 weight coordinates of $p(x)$ are partial derivatives of $q(x)$.
This identifies the space $G/P$ of type $\EE_6$ with the singular
locus of the cubic hypersurface $q(x)=0$.

\medskip

Let us define a symmetric bilinear form $p(x,y)$ on $V_1$
with values in $V_2$ by the formula $p(x+y)=p(x)+2p(x,y)+p(y)$.
Then $\exp(x+y)=\exp(x)\exp(y)$ implies that 
\begin{equation}
2p(x,y)=x\cdot y \label{a1}
\end{equation}
is the product of $x\in V_1$ and $y\in V_1$ in the 
commutative $k$-algebra $V$. 

We have a decomposition of 
$S^2(V_1)$ as the direct sum of $V_2$ and the representation
with highest weight $2\omega'$ (cf. the proof of Lemma \ref{p}).
In the notation of \cite{B}
the representation $V_2$ is irreducible with highest weight $\omega_1$,
in particular, it is minuscule. Thus
the eigenspaces for the action of the maximal torus $H'=H\cap G'$
are 1-dimensional, so that 
on $V_2$, in the same way as on $V_1$, 
we have weight coordinates well defined
up to a multiplicative constant. 
The coordinates $p_\lambda(x,y)$ of $p(x,y)$ are symmetric bilinear
forms of degree 2 with values in $k$. We can write 
\begin{equation}
p_\lambda(x,y)=\sum_{\lambda=\mu+\nu}p_{\mu\nu}x_\mu y_\nu, \label{ep}
\end{equation}
where $\mu$ and $\nu$ are weights of $V_1$, $p_{\mu\nu}\in k$,
and $x_\mu$ is a non-zero linear form on the weight $\mu$
subspace $(V_1)_\mu\subset V_1$ (and similarly for $y_\nu$ ). 
One checks that for $r=4,5,6,7$
the ranks of the quadratic forms
$p_\lambda(x)$ are 4, 6, 8, 10, respectively. 
If $r=7$ we associate to the 
cubic form 
$$q(x)=\sum_{\mu+\nu+\xi=0}q_{\mu\nu\xi}
x_\mu x_\nu x_\xi$$
the symmetric trilinear form
$$q(x,y,z)=\sum_{\mu+\nu+\xi=0}q_{\mu\nu\xi}
x_\mu y_\nu z_\xi.$$
In this case the weights of $V_2$ are the negatives of the weights of $V_1$.
Moreover, 
$$p_{-\mu}(x)=\frac{\partial q(x)}{\partial x_\mu},$$
so that 
\begin{equation}
3q(x,y,z)=\sum_\mu  p_{-\mu}(x,y) z_\mu, \quad
p_{-\mu}(x,y)=\sum_{-\mu=\nu+\xi}3q_{\mu\nu\xi}x_\nu y_\xi. \label{e12}
\end{equation}
For future reference we note that if $p_\lambda(x,y)=0$
for all $\lambda$, then $q(x,y,y)=0$. It follows from 
$\exp(x+y)=\exp(x)\exp(y)$ that 
\begin{equation}
3q(x,x,y)=p(x)\cdot y \label{a2}
\end{equation}
is the product of $p(x)\in V_2$ and $y\in V_1$ in the 
commutative $k$-algebra $V$.

\section{$G/P$ and blowing-up} \label{S4}

Let $\pi:(G/P)_a\to V_1$ be the restriction to $(G/P)_a$
of the natural projection $V=k\oplus V_1\oplus V_2\oplus V_3 \to V_1$.
We have $\exp(x)=(1,x,p(x),q(x))$ hence $\pi\circ\exp=\id$.
Here and in what follows we write our formulae for the case
$r=7$, with the convention that if $r<7$ the last
coordinate must be discarded.

We now describe the fibres of $\pi$.

\ble \label{fibre}
Let $g_t=(t,1,t^{-1},t^{-2})$, $t\in \ov k^*$. For 
$x\in V_1\otimes_k\ov k$ we have the following statements.

\noindent {\rm (a)} If $x\notin(G'/P')_a$, then
$\pi^{-1}(x)=\{g_t\cdot\exp(x)| t\in \ov k^*\}$.

\noindent {\rm (b)} If $x\in (G'/P')_a\setminus\{0\}$, then
$$\pi^{-1}(x)=\{(t,x,0,0)|t\in \ov k^*\}\cup
\{(0,x,2p_\lambda(x,u),3q(x,u,u))|u\in V_1\otimes_k\ov k\}.$$

\ele
{\it Proof} 
Recall that the torus $T$ is generated by the maximal torus 
$H\subset G$ and the scalar matrices $(t,t,t,t)$, $t\in \ov k^*$.
Let $h\in\h$ be an element of the Lie algebra of $H$ such that
$\beta (h)=0$ for all simple roots $\beta$ of $G$, $\beta\not=\alpha$, and
$\alpha (h)=1$. The 1-parameter subgroup $\G_m\subset H$
whose tangent vector 
at the identity is $h$, acts on $V$ as $(t^{m},t^{m-1},t^{m-2},t^{m-3})$,
where $m=\omega (h)$, and $\omega$ is the fundamental weight dual to $\alpha$.
Hence $g_t\in T$ for any $t\in \ov k^*$.

Every $\ov k$-point $y=(y_0,y_1,y_2,y_3)$ of the closed set $(G/P)_a$
satisfies the equations
\begin{equation}
y_0y_2=p(y_1), \quad y_0^2y_3=q(y_1), \label{e11}
\end{equation}
since these are satisfied on the affine cone over 
$\exp(V_1)$ which is dense in $(G/P)_a$.
Therefore, if $\pi$ sends a $\ov k$-point $y$ of $(G/P)_a$ to $x=y_1$, 
and $y_0\not=0$, 
we can write $y=g_t\cdot(1,x,p(x),q(x))=g_t\cdot\exp(x)$
for $t=y_0\in \ov k^*$. All such points are in $(G/P)_a$
since the action of $T$ preserves $(G/P)_a$, and $\exp(V_1)\subset(G/P)_a$.
If $y_0=0$ we see from (\ref{e11}) and Lemma \ref{p}
that $x\in (G'/P')_a$. This proves (a).

To prove (b) assume $x\in (G'/P')_a$, $x\not=0$. 
If $y_0\not=0$, then $y=(t,x,0,0)$, by (\ref{e11}).
 
We need some preparations for the case $y_0=0$.
Recall that $V_0$ is identified with $k$ by the choice 
of a highest weight vector $v\in V_0$, and $V_1$ is identified with
$\n$.
Consider $\g_1=\n^-$, the opposite nilpotent algebra of $\n$.
Any non-zero element $X\in \g_1$ sends $V_i$ to $V_{i-1}$ 
because of the grading. Hence we can write
$$\exp(Xt) (y_0,y_1,y_2,y_3)=(y_0+s(y_1,X)t+z_1t^2+z_2 t^3,
y_1+u_1t+u_2t^2, y_2+wt,y_3),$$
where $z_1,\,z_2 \in k$, $u_1,\,u_2\in V_1$, $w\in V_2$, and
$s(y_1,X)\in k$ is defined by 
$$s(y_1,X)v=X y_1 v=[X,y_1]v.$$
For any non-zero $y_1 \in \n\otimes_k\ov k=V_1\otimes_k\ov k$ 
one can find $X\in\g_1\otimes_k\ov k$ 
such that $s(y_1,X)=1$. Otherwise $\g_1 y_1 v=0$, and so
$y_1 v$ is a highest vector
of the $\g$-module $V\otimes_k\ov k$, which is not a multiple of $v$.
This contradicts the irreducibility of $V\otimes_k\ov k$. 
Let us fix such an element $X\in \g_1\otimes_k\ov k$. 

Now let $y_0=0$. Then
$$g_{t^{-1}}\, \exp(Xt)(0,y_1,y_2,y_3)=
(1+z_1 t+z_2 t^2,y_1+u_1t+u_2t^2,y_2t+wt^2,y_3t^2)$$
is a $\ov k[t]$-point of $(G/P)_a$, and hence
its coordinates satisfy (\ref{e11}) identically in $t$. 
Equating to $0$ the coefficient
at $t$ in the first equation in (\ref{e11}) we obtain
$y_2=2p(y_1,u)$, where $u=u_1$. 
Equating to $0$ the coefficient at $t^2$ 
in the second equation, and using that
$q(y_1,y_1,v)=0$ for all $v\in V_1$ (see (\ref{e12})) we obtain
$y_3=3q(y_1,u,u)$.

To complete the proof of (b) we need to show that for any $\ov k$-point
$x\in (G'/P')_a$ and any $u\in V_1\otimes_k\ov k$ the point
$(0,x,2p_\lambda(x,u),3q(x,u,u))$ is contained in $(G/P)_a$. We note that
$$(0,x,2p_\lambda(x,u),3q(x,u,u))=\exp(u)\cdot(0,x,0,0),$$
as immediately follows from (\ref{a1}) and (\ref{a2}).
Since $\exp(u)$ is in the unipotent group $N\subset G$ it is enough
to show that $(0,x,0,0)$ is in $(G/P)_a$. It is clear that
$(1,x,0,0)=\exp(x)$ is in $(G/P)_a$. Choosing $X\in \g_1\otimes_k\ov k$
as above such that $s(x,X)=-1$ we obtain $\exp(X)(1,x,0,0)=(0,x,0,0)$.
QED
\medskip

\bco \label{cc}
Let $U\subset (G/P)_a$ be the complement to the intersection of
$(G/P)_a$ with $(V_0\oplus V_1)\cup (V_2\oplus V_3)$.
The restriction of $\pi$ to $U$ is a morphism
$U\to V_1\setminus\{0\}$, which is the composition of a torsor
under the torus $\G_m=\{g_t|t\in \ov k^*\}$,
and the morphism inverse
to the blowing-up of $V_1\setminus\{0\}$ 
at $(G'/P')_a\setminus\{0\}$.
\eco
{\em Proof} The set $U$ is covered by the open subsets $U_0:y_0\not=0$,
and $U_\lambda:y_\lambda\not=0$, where the $y_\lambda$ are 
the weight coordinates in
$V_2$. Indeed, if $y_0=y_\lambda=0$ for all $\lambda$,
then we are in case (b) of Lemma \ref{fibre}, 
but $p_\lambda(x,u)=0$ for all $\lambda$ implies $q(x,u,u)=0$,
and such points are not in $U$. 
Each of these open subsets is $\G_m$-equivariantly isomorphic 
to the direct product of $\G_m$ and the closed subvariety 
of $(G/P)_a$ given by $y_i=1$ with trivial $\G_m$-action.
Gluing them together we obtain the quotient $\tilde U$. 

The equations (\ref{e11}) show that $\pi^{-1}(0)\cap U=\emptyset$,
thus $\pi$ projects $U$ to $V_1\setminus\{0\}$.
The action of $\G_m$ preserves the fibres, 
hence $\pi$ factors through a morphism $\tilde U\to V_1\setminus\{0\}$.
It is an isomorphism outside $(G'/P')_a$, whereas
the inverse image of $(G'/P')_a\setminus\{0\}$ is 
the projectivisation of the
normal bundle to $(G'/P')_a\setminus\{0\}$ in $V_1\setminus\{0\}$, 
by Lemma \ref{fibre} (b). It is not hard to prove
(and is well known to experts) that this implies that
$\tilde U$ is the blowing-up of $V_1\setminus\{0\}$ at
$(G'/P')_a\setminus\{0\}$. QED

\section{Del Pezzo surfaces} \label{DP}

For the geometry of exceptional curves on del Pezzo surfaces
the reader is referred to \cite{CF}, Ch. IV, see also \cite{FM}, Sect. 5.
Let $M_1,\ldots, M_r$, $4\leq r\leq 7$, be points in general
position in the projective plane $\P^2$, which says that no three 
points are on a line, and no six points are on a conic.
The blowing-up $X$ of $\P^2$ in $M_1,\ldots, M_r$ is called
a {\sl split del Pezzo surface} of degree $d=9-r$.
The surface $X$ contains exactly $d_r$ {\sl exceptional curves}, 
that is, smooth rational curves with self-intersection $-1$.
For $r\leq 6$ the exceptional curves on $X$ arise in
one of these ways:
the inverse images of the $M_i$; the proper transforms of the lines
through $M_i$ and $M_j$, $i\not=j$; the proper transforms of the
conics through five of the $M_i$. For $r=7$ one also has
the proper transforms of singular cubics passing through
all the 7 points with a double point at some $M_i$.
The intersection index
defines an integral bilinear form $(\,.\,)$ on $\Pic X$.
The opposite of the canonical class
$-K_X$ is an ample divisor, $(K_X^2)=d$. The Picard
group $\Pic \ov X=\Pic X$ is generated 
by the classes of exceptional curves (the complement to the
union of these curves is an open subset of $\A^2$).
The triple $(\Pic X, K_X, (\,.\,))$ coincides,
up to isomorphism, with the triple $(N_r, K_r, (\,.\,))$ defined as
$$N_r=\oplus_{i=0}^r\Z\ell_i, \ \  K_r=-3\ell_0+\sum_{i=1}^r\ell_i, 
\ \ (\ell_0^2)=1,
\ (\ell_i^2)=-1,\, i\geq 1, \, (\ell_i.\ell_j)=0,\, i\not=j,$$
see \cite[Thm. 23.9]{CF}. Moreover, 
the exceptional curves are identified with
the elements $\ell\in N_r$ such that $(\ell^2)=(\ell.K_r)=-1$,
called the {\sl exceptional classes} \cite[Thm. 23.8]{CF}.
By definition, a geometrically integral conic on $X$ is 
a smooth rational curve with self-intersection $0$. 
By Riemann--Roch theorem each conic
belongs to a 1-dimensional pencil of curves which are fibres of 
a morphism $X\to\P^1$, called a {\sl conic bundle}. We refer
to the fibres of such a morphism as {\sl conics}.
In particular, through every point of $X$
passes exactly one conic of a given pencil.
The classes of conic bundles can be
characterized by the properties $(c^2)=0$, $(c.K_r)=-2$.

Let $K_r^\perp$ be the orthogonal complement to $K_r$ in $N_r$.
The elements $\alpha\in K_r^\perp$ such that
$(\alpha^2)=-2$ form a root system $\RR$ 
in the vector space 
$K_r^\perp\otimes \R\simeq \R^r$ with the 
negative definite scalar product
$(\,.\,)$. In fact, $\RR$ is a root system
of rank $r$ in the series (\ref{series}). Moreover, the lattice
$K_r^\perp$ is generated by roots, so that $K_r^\perp\simeq Q(\RR)$.
For example, we can choose 
$$\beta_1=-\ell_1+\ell_2,\, \ldots,\, \beta_{r-1}=-\ell_{r-1}+\ell_r,\ 
\beta_r=-\ell_0+(\ell_1+\ell_2+\ell_3)$$ as a basis of simple roots of $\RR$.
The relation to our standard numeration (which follows \cite{B})
is this: $\alpha_r=\beta_{r-1}$, $\alpha_1=\beta_1$.

The Weyl group $\WW=\WW(\RR)$
generated by the reflections in the roots, 
is the automorphism group of the triple $(N_r,K_r,(\,.\,))$.
It operates transitively on the set of exceptional curves,
and also on the set of conic bundle classes, see, e.g.
\cite[Lemma 5.3]{FM}.
Let 
$$P(\RR)=\{n\in K_r^\perp\otimes \R|(n.m)\in\Z \ \text{for any}
\ m\in Q(\RR)\}$$
be the lattice dual to $Q(\RR)$; we have 
$Q(\RR)\subset P(\RR)$.
The image of the map 
$$N_r\to N_r\otimes\R=\R K_r\,\oplus\, (K_r^\perp\otimes\R)$$
is contained in the orthogonal direct sum
$\frac{1}{d}\Z K_r\oplus P(\RR)$ as a subgroup of index $d$.

\ble \label{omega}
Let $\alpha=\beta_{r-1}\in\RR$ be the simple root such that $(\RR,\alpha)$
is one of the pairs $(\ref{herm})$, and let $\omega\in P(\RR)$ be the
dual fundamental weight, $(\alpha.\omega)=-1$.

{\rm (i)} The exceptional classes in $N_r$ are 
$-\frac{1}{d}K_r+w\omega$, for all $w\in \WW$.

{\rm (ii)} Two distinct exceptional curves 
intersect in $X$ if and only if the corresponding
weights are not adjacent vertices of the convex hull 
$\Conv(\WW\omega)$.

{\rm (iii)} Let $\omega_1$ be the fundamental weight dual to the root 
$\beta_1$. The conic bundle classes in $N_r$ are 
$-\frac{2}{d}K_r+w\omega_1$, for all $w\in \WW$.
\ele

Note that since $\WW$ acts transitively on the set of bases, the choice of
a basis of simple roots is not important for the conclusion of this lemma.

\medskip

\noindent{\em Proof} (i) and (iii) 
The image of the exceptional class $\ell_r$ in $P(\RR)$ is 
the fundamental weight $\omega=\omega_{r-1}$, 
and the image of the conic bundle
class $\ell_0-\ell_1$ is the fundamental weight $\omega_1$.
The statement now follows from the transitivity of action of $\WW$
on these classes. Cf. \cite[Lemma 5.2]{FM}.

(ii) By the transitivity of $\WW$ on the 
exceptional classes it is enough to check 
this for the classes
$-\frac{1}{d}K_r+\omega$ and $-\frac{1}{d}K_r+x$,
where $x=w\omega$ for some $w\in \WW$. The intersection index 
\begin{equation}
(-\frac{1}{d}K_r+x.-\frac{1}{d}K_r+\omega)=\frac{1}{d}+(x.\omega)
\label{F}
\end{equation}
equals $-L(x)$ in the notation of the proof
of Proposition \ref{stable} (with the opposite sign of the scalar product).
In the simply laced case this proof shows that
$L(x)=1$ when $x=\omega$, $L(x)=0$ if $x$ is a
vertex of the convex hull 
$\Conv(\WW\omega)$ adjacent to $\omega$, and $L(x)<0$ for all other
$x\in \WW\omega$. QED

\medskip

We observe that for any conic bundle class
$x$ there exists a conic bundle class $y$ such that $(x.y)=1$.
Indeed, by the transitivity of $\WW$ on conic bundle classes
we can assume that $x=\ell_0-\ell_1$. For
$y=\ell_0-\ell_2$ we have $(x.y)=1$.

\section{Main theorem}\label{S5}

Let us recall our notation:

$(\RR,\alpha)$ be the pair in $(\ref{herm})$ such that 
$\RR$ has rank $r$;

$G$ is the simply
connected semisimple group with a split maximal torus $H$ and a 
maximal parabolic subgroup $P\supset H$, such that $(G,P)$
is defined by the pair $(\RR,\alpha)$;

$V$ is the fundamental representation of $G$ 
such that $P$ is the stabilizer of the line
spanned by a highest weight vector (this representation is
faithful);

$T\subset \GL(V)$ is the torus generated by the image of $H$ in $\GL(V)$,
and the scalar matrices;

$Y$ is the geometric quotient of $(G/P)_a^{sf}\subset(G/P)_a$
with respect to the natural left action of $T$;

the morphism $f : (G/P)_a^{sf}\to Y$ is a universal torsor (see Theorem \ref{univ}).

\noindent

Let $\Lambda\subset\hat H$ be the set of weights of $H$ in $V$, and let
$V_\lambda\subset V$ be the subspace of weight $\lambda$, so that
$V=\oplus_{\lambda\in\Lambda}V_\lambda$. In our case $\dim V_\lambda=1$
(since $V$ is minuscule, see Section \ref{S3}).
Let $\pi_\lambda:V\to V_\lambda$ be the natural projections,
and let $L_\lambda=\pi_\lambda^{-1}(0)$ be the weight coordinate
hyperplanes. For a subset $A\subset V$ we write $A^\times$
for the set of points of $A$ outside $\cup_{\lambda\in\Lambda}L_\lambda$.
For a subset $B\subset Y$ we write $B^\times$ for
$f(f^{-1}(B)^\times)$.

We now state our main theorem whose proof 
occupies the rest of the paper.

\bthe\label{main}
For $r=4$, $5$, $6$ or $7$
let $M_1,\ldots, M_r$ be points in general position in $\P^2$
(no three on a line, no six on a conic).
Let $X$ be the blowing-up of $\P^2$ in $M_1,\ldots, M_r$.
There exists an embedding $X\hookrightarrow Y$ such that
$X\setminus X^\times$ is the union of exceptional curves on $X$. 
For such an embedding $f^{-1}(X)\to X$ is a universal torsor.
\ethe

We write $S^n_\chi(V)$ for the weight $\chi\in \hat H$ subspace
of $S^n(V)$, and $S^n_\chi(V)^*$ for the dual
space of functions. Let $I(\T)\subset k[V]=S(V^*)$ be the ideal of $\T$.
We shall prove the following statement from which the main theorem
will follow:

\medskip

{\it There exists an embedding of a universal torsor
$\T$ over $X$ into $(G/P)^{sf}_a\subset V$ such that the restriction
of $f$ to $\T$ is the structure morphism $\T\to X$, and 
$f(\T^\times)$ is the complement to the union of exceptional curves on $X$.
Moreover, for $r<7$ the ideal $I(\T^\times)\subset k[V^\times]$ 
is generated by the graded 
components of degree $2$ and weight $w\omega_1$, for all $w\in \WW$.}

\medskip

The last statement will be used in the case $r=7$, and can be ignored
by the reader interested in the cases $r=5$ and $r=6$ only. Recall that $\omega_1$
is the highest weight of a non-trivial irreducible $\g$-module of least dimension.

\medskip

\noindent{\em Proof} The proof is by induction on $r$ starting from
$r=4$. In this case $Y$ is a del Pezzo surface of degree 5,
$G/P$ is the Grassmannian variety $G(3,5)\simeq G(2,5)$,
and $G(3,5)^{sf}=G(3,5)^s$ is a universal torsor over $Y$
(see \cite{S1} or \cite[Lemma 3.1.6]{S}).
It is well known that the ideal of $G(3,5)_a\subset V$ is generated by
the (quadratic) Pl\"ucker relations, and it is easy to see that
their weights are of the form $w\omega_1$, 
so that our statement is true in this case.

Suppose we know the statement for $r-1\geq 4$. This means that
we are given the following data: 

$(\RR',\alpha')$ is the `previous' pair to $(\RR,\alpha)$ in 
$(\ref{herm})$;

$\WW'=\WW(\RR')$ is the Weyl group;

$G'$ and $P'$ are defined by $(\RR',\alpha')$,
so that $(G'/P')_a\subset V_1$ (see Section \ref{S3});

$H'=H\cap G'$, so that 
$\RR'$ is the root system of $G'$ with respect to $H'$;

$T'\subset \GL(V_1)$ is the torus generated by the image
of $H'$ in $\GL(V_1)$ and the scalars ($T'$ is also the image
of $H$ in $\GL(V_1)$);

$x_\mu$ is a non-zero linear form on the weight $\mu$ subspace of $V_1$;

$Y'$ is the quotient of $(G'/P')_a^{sf}$ by $T'$;

$f' : (G'/P')_a^{sf}\to Y'$ is a universal torsor;

$X'$ is the blowing-up of $\P^2$ in $M_1,\ldots, M_{r-1}$
(it is a del Pezzo surface of degree $d'=8-r$);

an embedding $X'\hookrightarrow Y'$ satisfying the conditions
of the theorem, in particular,

$\T'=f'^{-1}(X')\to X'$ is a universal torsor.

The general position assumption implies that $M_r$
does not belong to the exceptional curves of $X'$.
Thus, by Hilbert's theorem 90, we can find a $k$-point
$x_0\in \T'^\times$ such that $f'(x_0)=M_{r}$.

\smallskip

Let $\tau:\hat T'\to\Pic X'$ be the map
defined in Section \ref{S1}; up to sign $\tau$ coincides
with the type of the torsor $f':\T'\to X'$ (Lemma \ref{tau}).
Since the torsor $f':\T'\to X'$ is universal,
$\tau$ is an isomorphism of $\hat T'=K[\T'_K]^*/K^*$ and $\Pic X'$
as abelian groups. To account for the duality between 
vectors and linear forms on $V_1$
we identify these groups by the isomorphism $-\tau$.
Recall that the Weyl group $\WW'$ acts on $\hat T'$
via the normalizer of $H'$ in $G'$, permuting the weights of $V_1$. 
By induction assumption $-\tau$
sends these weights bijectively onto the exceptional classes
in $\Pic X'$. If we transport the action of $\WW'$ from $\hat T'$
to $\Pic X'$ using $-\tau$, then the action of $\WW'$ so obtained
preserves the intersection index of exceptional curves, see (\ref{F}).
Thus $-\tau$ is a homomorphism of $\WW'$-modules,
where $\WW'$ acts on $\Pic X'$ as
the automorphism group of the triple $(N_{r-1},K_{r-1},(\,.\,))$.
In particular, $-\tau$ identifies the $\WW'$-(co)invariants
on both sides (isomorphic to $\Z$). 
This implies that if $\chi$ is a weight of $T'$ in 
$S^n(V_1)$, then the restriction of $\chi$
to the scalar matrices $\G_m\subset T'$ coincides with the intersection index
of $-\tau(\chi)$ with $-K_{X'}$, that is, 
\begin{equation}
(\tau(\chi).K_{X'})=n \label{6.2}
\end{equation}
(the sign is uniquely determined by the fact that
effective divisors intersect positively with $-K_{X'}$).
The isomorphism $-\tau$ also identifies the quotients by
the $\WW'$-invariants, that is, $P(\RR')$ and $\hat H'$.
We fix these identifications from now on.

For $\phi(x)\in S^n_\chi(V_1)^*$, $\chi\in\hat T'$, we let
$C_\phi\subset X'$ be the image of the intersection of $\T'$
with the $T'$-invariant hypersurface $\phi(x)=0$.
If $C_\phi\not=X'$, then the class $[C_\phi]$ in $\Pic X'$
is $-\tau(\chi)$, and (\ref{6.2}) can be written as
\begin{equation}
([C_\phi].-K_{X'})=n. \label{6.2a}
\end{equation}
We have (see the end of Section \ref{S1} for the first equality)
\begin{equation}
\H^0(X',\O_{-\chi})=k[\T']_{-\chi}=S^n_\chi(V_1)^*/I(\T')\cap S^n_\chi(V_1)^*.
\label{ideal}
\end{equation}
Apart from the weights of $V_1$ which correspond to exceptional
curves, the following two cases will be particularly relevant.
For $n=2$ let $\lambda$ be a weight of $T'$ in $V_2$. 
The restriction of $\lambda$ to $H'$ is 
$w\omega_1\in \hat H'=P(\RR')$,
where $w\in \WW'$ (see the end of Section \ref{S3}).
If $\phi\in S^2_\lambda(V_1)^*$ is such that $C_\phi\not=X'$, then
by (\ref{6.2a}) we see
that $[C_\phi]=-\frac{2}{d'}K_{X'}+w\omega_1$, so
$C_\phi$ is a conic on $X'$ by Lemma \ref{omega} (iii). The
Riemann--Roch theorem implies that $\dim \H^0(X',\O_{-\lambda})=2$,
where $\O_{-\lambda}=\O(C_\phi)$ is the invertible sheaf associated to $C_\phi$.
Thus $I(\T')\cap S^2_\lambda(V_1)^*$ has codimension 2 in $S^2_\lambda(V_1)^*$.
Note that by Lemma \ref{p} we have $p_\lambda(x)\in I(\T')\cap S^2_\lambda(V_1)^*$.

For $r=7$ and $n=3$ the space $V_3$ is a trivial $1$-dimensional 
representation of $G'$, hence of weight $0\in\hat H'$.
Thus for $\phi\in S^3_{0}(V_1)^*$
we have $[C_\phi]=-K_{X'}$, by (\ref{6.2a}). If $C_\phi\not=X'$,
then $C_\phi$ is a plane section of the cubic surface $X'\subset\P^3$. 
The vector space $\H^0(X',\O(C_\phi))=\H^0(X',\O(-K_{X'}))$
has dimension 4, thus 
$I_0=I(\T')\cap S^3_{0}(V_1)^*$ has codimension 4
in $S^3_{0}(V_1)^*$. It is clear that $q(x)\in I_0$, see, e.g.,
(\ref{a2}).

The following proposition is a crucial technical step in the
proof of our main theorem.

\bpr
There exists a non-empty open subset $\Omega(x_0)\subset (G'/P')_a^\times$
such that for any $y_0\in \Omega(x_0)$ 
we have $p_\lambda(x_0^{-1}y_0x)\notin I(\T')\cap S^2_\lambda(V_1)^*$
for all weights $\lambda$ of $V_2$,
and $q(x_0^{-1}y_0x)\notin I_0$ if $r=7$.
\epr
{\em Proof} We begin with pointing out the following useful fact.
Let ${\rm Ver}_\lambda$ be the composition of 
the second Veronese embedding
$V_1\to S^2(V_1)$ with the projection of 
$S^2(V_1)$ to its direct summand $S^2_\lambda(V_1)$.
By Lemma \ref{p}, $p_\lambda(x)=0$ is the only quadratic equation
of $G'/P'$ of weight $\lambda$,
thus ${\rm Ver}_\lambda((G'/P')_a)$ spans a codimension 1 subspace of
$S^2_\lambda(V_1)$, namely, the zero set of the linear form
$p_\lambda(x)\in S^2_\lambda(V_1)^*$.

Next, we claim that the quadratic forms
$p_\lambda(x_0^{-1}y_0x)$, $y_0\in (G'/P')_a^\times$, span
a codimension 1 subspace of $S^2_\lambda(V_1)^*$.
Using (\ref{ep}) we write
$$p_\lambda(x_0^{-1}y_0x)=\sum_{\lambda=\mu+\nu}
p_{\mu\nu}\frac{y_{0\mu} y_{0\nu}}{x_{0\mu} x_{0\nu}}x_\mu x_\nu.$$
Suppose that for some coefficients $c_{\mu\nu}$ we have a linear
relation
$$\sum_{\lambda=\mu+\nu}c_{\mu\nu}p_{\mu\nu}
\frac{y_{0\mu} y_{0\nu}}{x_{0\mu} x_{0\nu}}=0.$$
This can be read as a relation with coefficients
$c_{\mu\nu}p_{\mu\nu}x_{0\mu}^{-1} x_{0\nu}^{-1}$ satisfied
by all the vectors $(y_{0\mu} y_{0\nu})$, where
$y_0\in (G'/P')_a^\times$ and $\mu+\nu=\lambda$. 
The set of these vectors is
precisely ${\rm Ver}_\lambda((G'/P')_a^\times)$.
The linear span of ${\rm Ver}_\lambda((G'/P')_a^\times)$
is the same as the linear span of ${\rm Ver}_\lambda((G'/P')_a)$.
By the argument in the beginning of the proof,
up to a multiplicative constant there is
only one linear relation satisfied by the elements of
${\rm Ver}_\lambda((G'/P')_a)$, namely the one with coefficients $p_{\mu\nu}$.
Therefore, $c_{\mu\nu}=x_{0\mu} x_{0\nu}$ is uniquely determined up to
a multiplicative constant. This proves our claim.
Note that the linear span under discussion is thus
the space of forms vanishing at $x_0$.

It follows that the set of $\ov k$-points $y\in (G'/P')_a^\times$ such that
$p_\lambda(x_0^{-1}y x)$ belongs to the codimension 2 subspace 
$I(\T')\cap S^2_\lambda(V_1)^*$,
is a proper closed subset of $(G'/P')_a^\times$. 
For $r<7$ we define $\Omega(x_0)$ as the complement
to the union of these closed subsets for all weights $\lambda$ of $V_2$.

Until the rest of the proof we let $r=7$. 
Let ${\rm Ver}^3_0:V_1\to S^3_{0}(V_1)$ be
the composition of the natural map
$V_1\to S^3(V_1)$ with the projection $S^3(V_1)\to S^3_{0}(V_1)$.
The map ${\rm Ver}^3_0$ sends $x=(x_\mu)$ to the vector 
$(x_\mu x_\nu x_\xi)$, for all $\mu,\,\nu,\,\xi$ such that $\mu+\nu+\xi=0$.
If we write the invariant cubic form (defined up to a scalar
multiple) as 
$$q(x)=\sum_{\mu+\nu+\xi=0} q_{\mu\nu\xi}x_\mu x_\nu x_\xi,$$
then it is well known that all the coefficients
$q_{\mu\nu\xi}$ are non-zero (see, e.g., \cite{F}). Recall that
the singular locus of the cubic hypersurface $q(x)=0$ is $(G'/P')_a$. 

Let $L_{x_0}\subset S^3_{0}(V_1)^*$ be the subspace of forms
vanishing at $x_0$ together with all their (first order) partial derivatives.
We claim that $L_{x_0}$ coincides with the linear span of the forms $q(x_0^{-1}y_0u)$,
where $y_0$ ranges over $(G'/P')_a$.

Let us prove this claim.
The partial derivatives of $q(x)$ vanish on 
$(G'/P')_a$, hence $q(x_0^{-1}y_0u)\in L_{x_0}$ for any $y_0\in (G'/P')_a$.
Thus the linear span of the forms $q(x_0^{-1}y_0u)$,
where $y_0\in(G'/P')_a$, is contained in $L_{x_0}$.
We now prove that these spaces have the same dimension.

Let $f(x)=\sum_{\mu+\nu+\xi=0} f_{\mu\nu\xi}x_\mu x_\nu x_\xi$ 
be a form in $L_{x_0}$.
The partial derivative with respect to $x_\xi$ is
$3\sum_{\mu+\nu=-\xi} f_{\mu\nu\xi}x_\mu x_\nu$. It vanishes at $x_0\in V_1^\times$
if and only if 
$$x_\xi\sum_{\mu+\nu=-\xi} f_{\mu\nu\xi}x_\mu x_\nu=
\sum_{\mu+\nu=-\xi}q_{\mu\nu\xi}^{-1}f_{\mu\nu\xi} \cdot 
q_{\mu\nu\xi}x_\mu x_\nu x_\xi$$ 
does. Hence $(q_{\mu\nu\xi}^{-1}L_{x_0})^\perp$ is spanned by the 27 vectors
$(q_{\mu\nu\xi}x_{0\mu} x_{0\nu} x_{0\xi})$, 
where $\xi$ is fixed, and $\mu,\,\nu$ are
arbitrary. Since the coordinates of $x_0$ are non-zero,
this space has the same dimension as the space 
$M\subset S^3_{0}(V_1)$ 
spanned by the 27 vectors $(q_{\mu\nu\xi})$, where $\xi$ is fixed, 
and $\mu,\,\nu$ are arbitrary weights satisfying $\mu+\nu+\xi=0$. 
The fact that the ideal of $(G'/P')_a$
is generated by the partial derivatives of $q(x)$, implies that
$M^\perp$ is the linear span of 
${\rm Ver}^3_0((G'/P')_a)$. We conclude that $\dim L_{x_0}$
equals the dimension of this linear span.
Since all the coefficients $q_{\mu\nu\xi}$ are non-zero, the 
forms $q(x_0^{-1}y_0u)$,
where $y_0\in(G'/P')_a$, span the space of the same dimension.
This proves our claim.

Let us complete the proof of the proposition in the case $r=7$.
A cubic form $f\in S^3_{0}(V_1)^*$ is in $L_{x_0}$
if and only if $f(x)=0$
is singular at $x_0\in V_1^\times$. This is the case if and only if
the corresponding hyperplane $H_f\subset S^3_{0}(V_1)$ contains
the tangent space $\Phi$
to ${\rm Ver}^3_0(V_1)$ at the point $m={\rm Ver}^3_0(x_0)$.
We have a commutative diagram (cf. (\ref{ideal}) and (\ref{end}))
$$\begin{array}{ccccc}
X'&\leftarrow&\T'&\hookrightarrow &V_1 \\
\downarrow&&\downarrow&&\downarrow \\
\P(\H^0(X',\O(-K_{X'}))^*)&\leftarrow&
\H^0(X',\O(-K_{X'}))^*\setminus\{0\}&\hookrightarrow& 
S^3_{0}(V_1)
\end{array}$$
where the left hand vertical map is the anticanonical embedding of $X'$,
and the other two are ${\rm Ver}^3_0$. 
The image of $\T'$ in the 4-dimensional vector space
$$\H^0(X',\O(-K_{X'}))^*=(k[\T']\cap S^3_0(V_1)^*)^*=
(S^3_0(V_1)^*/I_0)^*\simeq \A^4\,\subset S^3_0(V_1)$$ 
is the affine cone $X'_a$ (without $0$) over the cubic
surface $X'\subset\P^3$. 

By induction assumption $I(\T'^\times)$ is 
generated by its graded components $I_\lambda$
of degree $2$ and weight $\lambda$, for all weights $\lambda$ 
of $V_2$. The weights of $V_1$ are the negatives of the weights of
$V_2$, so that $x_{-\lambda}I_{\lambda}$ has degree 3 and weight 0.
Since the coordinates $x_{-\lambda}$ are invertible on $\T'^\times$,
the ideal $I(\T'^\times)$ is generated by its graded component
of degree 3 and weight 0. Hence locally in the neighbourhood 
$\T'^\times$ of $x_0$ 
the ideal $I(\T')$ is generated by $I_0$, that is,
by the equations of $\A^4$ in $S^3_{0}(V_1)$.

This implies that the tangent space $T_{X'_a,m}\subset\A^4$
is $\Phi\cap\A^4$. Thus for any $f$ in a dense open subset 
of $L_{x_0}$ we have $H_f\cap\A^4=T_{X'_a,m}$.
Since $X'\subset\P^3$ is a smooth cubic surface,
$X'_a\setminus T_{X'_a,m}$ is dense and open in $X'_a$.
Therefore, for the general $f\in L_{x_0}$ we have 
$X'_a\cap H_f\not= X'_a$, so that $f\notin I_0$.
Now the above claim implies the statement of proposition. QED

\bco \label{lemma1}
For any $k$-point $y_0\in \Omega(x_0)$ and any weight $\lambda$ of $V_2$
the closed subset of $\T'$ given by $p_\lambda(x_0^{-1}y_0 x)=0$
is the preimage $f'^{-1}(C_\lambda)$
of a geometrically integral $k$-conic $C_\lambda\subset X'$ 
passing through $M_{r}$.
For $r=7$ the closed subset of $\T'$ given by 
$q(x_0^{-1}y_0x)=0$, for any $y_0\in \Omega(x_0)$, is the preimage $f'^{-1}(Q)$
of a geometrically integral cubic $k$-curve $Q$ 
with a double point at $M_r$ (the intersection of the cubic surface
$X'$ with its tangent plane at $M_r$).
\eco
{\em Proof} To check that $M_{r}\in C_\lambda$ set $x=x_0$, then 
$p_\lambda(x_0^{-1}y_0x)=p_\lambda(y_0)=0$ by Lemma \ref{p} 
since $y_0\in(G'/P')_a$. 
If the conic $C_\lambda$ is not geometrically integral, then
its components must have intersection index $1$ with $-K_{X'}$,
so there are two of them. It is well known that a curve on $X'$ 
has such a
property if and only if it is an exceptional curve.
However, $M_r$ does not belong to
the exceptional curves of $X'$. Thus $C_\lambda$ is geometrically integral.
 
If $r=7$, by substituting $x=x_0$ one shows as before that $Q$
contains $M_7$ (the cubic form $q$ vanishes on $G'/P'$).
Since the $p_\lambda(x)$ are partial derivatives of $q(x)$, 
and $M_7\in C_\lambda$, we see that $Q$ has a double point at $M_7$. 
If $Q$ is not geometrically integral, then it is the union of 
a geometrically integral conic and an exceptional curve, 
or the union of three
exceptional curves. In each of these cases the singular point
$M_7\subset Q$ will have to lie
on an exceptional curve, and this is a contradiction. QED

\bco \label{cccc}
For any $y_0\in \Omega(x_0)$ the scheme-theoretic
intersection of $x_0^{-1}y_0\T'$ and $(G'/P')_a$ is the orbit $T'y_0$.
\eco 
{\em Proof} By Lemma \ref{p} the ideal of $(G'/P')_a$ 
is generated by $p_\lambda(x)$, for all weights $\lambda$ of $V_2$.
As was remarked at the end of Section \ref{DP}, there exist
weights $\lambda$ and $\nu$ such that the intersection index of 
$C_\lambda$ and $C_\nu$ on $X'$ is $1$, that is, $M_r$ is the 
scheme-theoretic intersection $C_\lambda \cap C_\nu$.
Thus the orbit $T'y_0$ is the closed subscheme 
of $x_0^{-1}y_0\T'$ given by $p_\lambda(x)=p_\nu(x)=0$,
and our statement follows. QED

\medskip

Let $\sigma : X=\Bl_{M_r}(X')\to X'$ be the morphism 
inverse to the blowing-up of $M_r$.
Then $\sigma$ induces an isomorphism of $X\setminus \sigma^{-1}(M_r)$
with $X'\setminus M_r$, and $\sigma^{-1}(M_r)\cong \P^1$.
The proper transform of a curve $D\subset X'$ 
is defined as the closure of $\sigma^{-1}(D\setminus M_{r})$ in $X$.
The comparison of intersection indices on $X'$ and $X$ shows
that the proper transforms of the conics $C_\lambda$ and the singular
cubic $Q$ (for $r=7$)
are exceptional curves on $X$. 
By comparing the numbers we see that these curves together with
$\sigma^{-1}(M_r)$ and the inverse images of the exceptional curves on $X'$
give the full set of exceptional curves on $X$.

\medskip

\noindent{\em End of proof of Theorem $\ref{main}$.}
Consider the open set $U\subset (G/P)_a$ 
and the morphism $\pi:U\to V_1\setminus\{0\}$, see Corollary \ref{cc}.
Choose any $y_0\in \Omega(x_0)$, and define $\T\subset U$ 
 as the `proper
transform' of $x_0^{-1}y_0\T'$ with respect to $\pi$.
Explicitly, $\T\subset U$ is defined as the Zariski closure of 
$$\pi^{-1}(x_0^{-1}y_0\T'\setminus (G'/P')_a)
=\pi^{-1}(x_0^{-1}y_0\T'\setminus T'y_0),$$
where the equality is due to Corollary \ref{cccc}. 
The torus $T'$ acts on $\T'$, and
$\pi$ is $T'$-equivariant, hence $T'$ acts on $\T$.
But $\G_m=\{g_t\}$ (see Lemma \ref{fibre}) also
acts on $\T$. The torus $T$ is generated by $T'$ and $\G_m=\{g_t\}$,
so that $T$ acts on $\T$.

Corollaries \ref{cc} and \ref{cccc} imply that 
the restriction of $\pi$ to $\T$
is the composition of a torsor under $\G_m=\{g_t\}$
and the morphism
$\Bl_{y_0T'}(x_0^{-1}y_0\T')\to x_0^{-1}y_0\T'$ 
inverse to the blowing-up of the orbit $T' y_0$ in $x_0^{-1}y_0\T'$.
The blowing-up of $T' y_0$ in $x_0^{-1}y_0\T'$ is naturally isomorphic to
the pullback $\T'\times_{X'}X$ of the torsor $\T'\to X'$ to $X$.
This can be summarized in the following commutative diagram:
\begin{equation}
\begin{array}{ccccc}
\T&\hfl{}{}{6mm} &\T'\times_{X'}X &\hfl{}{}{6mm} &X\\
&&\vfl{}{}{5mm}&&\vfl{\sigma}{}{5mm}\\
&&\T'&\hfl{}{}{6mm} &X'
\end{array}
\label{last}
\end{equation}
where the horizontal arrows are torsors under tori,
and the vertical arrows are contractions.
The composed morphism $f:\T\to X$ is a composition
of two torsors under tori, and hence is an affine morphism
whose fibres are orbits of $T$. Therefore
$\T$ is an $X$-torsor under $T$, by Lemma \ref{torsor}.
We obtain a $T$-equivariant embedding $\T\hookrightarrow (G/P)_a$.

For $r<7$ we note that $I(\T^\times)\subset k[V^\times]$ 
is generated by $I(x_0^{-1}y_0\T'^\times)$ 
and the equations of $(G/P)_a$, moreover, for each weight
$w\omega_1$, $w\in \WW$, there is exactly one quadratic equation,
by Lemma \ref{p}. The restriction of $\omega_1\in \hat H=P(\RR)$
to $H'$ is again the weight $\omega_1\in \hat H'=P(\RR')$.
By induction assumption $I(\T'^\times)$ is generated by 
its graded components of degree 2 of such weights, 
hence the same is true for $I(\T^\times)$.

It remains to prove that $\T\subset (G/P)_a^{sf}$, and that the torsor
$f:\T\to X$ is universal.
The action of $T$ on $\T$ is free, so let us show that every point of
$\T$ is stable. We claim that $f$ sends 
the weight hyperplane sections of $\T$ to the exceptional curves on $X$.
By the results of Section \ref{S4} this follows from induction assumption
for the weights of $V_1$,
and from Corollary \ref{lemma1} for the weights of $V_2\oplus V_3$. 
Corollary \ref{cccc} implies that the highest weight hyperplane 
$x_\omega=0$ corresponds
to $\sigma^{-1}(M_r)$. By Lemma \ref{omega} (ii)
the set of exceptional curves of $X$
is identified with the set $\WW \omega$ in such a way
that two distinct exceptional curves 
intersect in $X$ if and only if the corresponding
weights are not adjacent vertices of the convex hull 
$\Conv(\WW\omega)$. Now Proposition \ref{stable}
implies that $\T\subset (G/P)^{sf}_a$. We thus obtain
an embedding $X\hookrightarrow Y$.

The pull-back of the torsor $(G/P)^{sf}_a\to Y$
to $X$ gives rise to
the following commutative diagram, where the horizontal
arrows represent the types of corresponding torsors:
$$
\begin{array}{ccc}
\hat T  & \hfl{\sim}{}{6mm} & \Pic Y\\
||&&\downarrow\\
\hat T& \hfl{}{}{6mm} & \Pic X
\end{array}
$$
The upper horizontal arrow is an isomorphism
since the torsor $(G/P)^{sf}_a\to Y$ is universal,
by Theorem \ref{univ}. Since
the exceptional curves on $X$ are cut by divisors on $Y$, 
the restriction map $\Pic Y\to\Pic X$
is surjective. However, the ranks of $\Pic Y$ and $\Pic X$
are equal, so this map is an isomorphism.
Now it follows from the diagram that
the type of the torsor $f:\T\to X$
is an isomorphism, so that this torsor is universal as well.
The theorem is proved. QED

\noindent Department of Mathematics, University of California,
Berkeley, CA, 94720-3840 USA
\medskip

\noindent serganov@math.berkeley.edu

\bigskip

\noindent Department of Mathematics, South Kensington Campus, 
Imperial College London, SW7 2BZ England, U.K.

\smallskip

\noindent Institute for the Information Transmission Problems, 
Russian Academy of Sciences, 19 Bolshoi Karetnyi, 
Moscow, 127994 Russia
\medskip

\noindent a.skorobogatov@imperial.ac.uk


\begin{thebibliography}{99}

\bibitem{BP} V.V. Batyrev and O.N. Popov. The Cox ring of a 
del Pezzo surface. In: {\it Arithmetic of 
higher-dimensional algebraic varieties} 
(Palo Alto, 2002), Progr. Math. {\bf 226}
Birkh\"auser, 2004, 85--103.

\bibitem{B} N. Bourbaki. {\it Groupes et alg\`ebres de Lie.} 
Chapitres IV-VIII. Masson, Paris, 1975, 1981.

\bibitem{CS}
J-L. Colliot-Th\'el\`ene et J-J. Sansuc. La descente sur
les vari\'et\'es rationnelles, II. {\it Duke Math. J.} {\bf 54}
(1987) 375--492.

\bibitem{Da}
R. Dabrowski. On normality of the closure of a generic torus orbit in 
$G/P$. {\it Pacific J. Math.} {\bf 172} (1996) 321--330.

\bibitem{D} U. Derenthal. Universal torsors of Del Pezzo surfaces 
and homogeneous spaces. {\it Adv. Math.}, to appear. arXiv:math.AG/0604195

\bibitem{D1} U. Derenthal. On the Cox ring of Del Pezzo surfaces.
arXiv:math.AG/0603111

\bibitem{Do}
I.V. Dolgachev. {\it Lectures on invariant theory.} 
Cambridge University Press, 2003.

\bibitem{FH}
H. Flaschka and L. Haine. Torus orbits in $G/P$. {\it Pacific J. Math.} 
{\bf 149} (1991) 251--292.

\bibitem{F}
J. Faulkner. Generalized quadrangles and cubic forms. {\it Comm. Algebra}
{\bf 29} (2001) 4641--4653.

\bibitem{FM}
R. Friedman and J.W. Morgan. Exceptional groups and del Pezzo surfaces. 
In: Symposium in Honor of C. H. Clemens (Salt Lake City, UT, 2000), 
101--116, {\it Contemp. Math.} {\bf 312} Amer. Math. Soc., Providence, RI, 2002.

\bibitem{GS}
I.M. Gelfand and V.V. Serganova. Combinatorial geometries 
and the strata of a torus on homogeneous compact manifolds. 
(Russian) {\it Uspekhi Mat. Nauk} {\bf 42} (1987) 107--134.

\bibitem{H}
R. Hartshorne. {\it Algebraic geometry.} Springer-Verlag, 1977.

\bibitem{He}
S. Helgason. {\it Differential geometry, Lie groups, and symmetric spaces.}
Graduate Studies in Math. {\bf 34} Amer. Math. Soc. 1978, 2001.

%\bibitem{KKV}
%F. Knop, H. Kraft and T. Vust. The Picard group of a $G$-variety. 
%In: {\it Algebraische Transformationsgruppen und Invariantentheorie}, 
%DMV Sem. {\bf 13} Birkh\"auser, Basel, 1989, 77--87.

\bibitem{LT} G. Lancaster and J. Towber.
Representation-functors and flag-algebras for the classical groups. I. {\it J. Algebra}
{\bf 59} (1979) 16--38.

\bibitem{Leung}
N.C. Leung. ADE-bundles over rational surfaces, configuration of lines and rulings. math.AG/0009192

\bibitem{CF} Yu.I. Manin. {\it Cubic forms.} 2nd. ed. North-Holland, 1986.

\bibitem{M} L. Manivel. Configurations of lines and models of Lie algebras.
{\it J. Algebra} {\bf 304} (2006) 457--486.

\bibitem{GIT}
D. Mumford, J. Fogarty, and F. Kirwan. {\it Geometric invariant theory.}
3rd enlarged edition. Springer-Verlag, 1994.

\bibitem{VO} 
A.L. Onishchik and E.B. Vinberg. {\it Lie groups and algebraic groups.}
Springer-Verlag, 1990.

\bibitem{VP} V.L. Popov. Picard groups of homogeneous spaces of linear 
algebraic groups and one-dimensional homogeneous vector fiberings. 
(Russian) {\it Izv. Akad. Nauk SSSR} Ser. Mat. {\bf 38} (1974) 294--322.

\bibitem{P} O.N. Popov. {\it Del Pezzo surfaces and algebraic groups.}
Diplomarbeit, Universit\"at T\"ubingen, 2001.

\bibitem{S}
A.N. Skorobogatov. {\it Torsors and rational points.} 
Cambridge tracts in mathematics {\bf 144},
Cambridge University Press, 2001.

\bibitem{S1}
A.N. Skorobogatov. On a theorem of Enriques--Swinnerton-Dyer. {\it Ann.
Fac. Sci. Toulouse} {\bf 2} (1993) 429--440.


\end{thebibliography}
\end{document}